\documentclass{siamart190516}

\usepackage{hyperref}
\usepackage{epstopdf}
\usepackage{graphicx}
\usepackage{graphicx}
\usepackage{amsmath}
\usepackage{amssymb,bm}

\usepackage{algorithm}
\usepackage{algorithmic}

\graphicspath{ {./images/} }

\usepackage{comment}

\newcommand{\cc}[1]{\mathcal{#1}}
\newcommand{\bb}[1]{\mathbb{#1}}

\newcommand{\bx}{\bm x}
\newcommand{\by}{\bm y}

\newcommand{\bA}{\bm A}

\newcommand{\cA}{\cc A}
\newcommand{\cX}{\cc X}
\newcommand{\cG}{\cc G}
\newcommand{\cF}{\cc F}
\newcommand{\cB}{\cc B}
\newcommand{\cY}{\cc Y}

\newcommand{\cU}{\cc U}
\newcommand{\cV}{\cc V}
\newcommand{\cS}{\cc S}

\newcommand{\RR}{\mathbb R}

\usepackage{textcomp}
\usepackage{adjustbox}

\begin{document}

\title{TT-LSQR for tensor least squares problems and application to data mining
\thanks{Version of \today. The authors are members of the INdAM 
Research Group GNCS that partially supported this  work.} }{}
\date{\today}
\author{ Lorenzo Piccinini\thanks{Dipartimento di Matematica, 
        Alma Mater Studiorum - Universit\`a di Bologna, Piazza di Porta San Donato 5,
        40126 Bologna, Italia. Email:
{\tt  lorenzo.piccinini12@unibo.it}}\and Valeria Simoncini\thanks{Dipartimento di Matematica, AM$^2$,
        Alma Mater Studiorum - Universit\`a di Bologna, Piazza di Porta San Donato 5,    40126 Bologna, Italy, and IMATI-CNR, Via Ferrata 5/A, Pavia, Italy. Email:
{\tt valeria.simoncini@unibo.it}}}

\maketitle
\begin{center}
    {\it This paper is dedicated to Michela Redivo-Zaglia and  Hassane Sadok}
\end{center}

\begin{abstract}
We are interested in the numerical solution of the tensor least squares problem 
\begin{equation*}
    \min_{\cc{X}} \| \cc{F} - \sum_{i =1}^{\ell} \cc{X} \times_1 A_1^{(i)} \times_2 A_2^{(i)} \cdots \times_d A_d^{(i)} \|_F,
\end{equation*}
where $\cc{X}\in\bb{R}^{m_1 \times m_2 \times \cdots \times m_d}$, $\cc{F}\in\bb{R}^{n_1\times n_2 \times \cdots \times n_d}$ are tensors with $d$ dimensions, and
the coefficients $A_j^{(i)}$ are tall matrices of conforming dimensions. We first
describe a tensor implementation of the classical 
LSQR method by Paige and Saunders, using the tensor-train
representation as key ingredient. We also show how to incorporate
sketching to lower the computational cost of dealing
with the tall matrices $A_j^{(i)}$.
We then use this methodology
to address a problem in information retrieval, the 
classification of a new query document among already
categorized documents, according to given keywords.
\end{abstract}



\section{Introduction}\label{sec:intro}

We are interested in the numerical solution of the multiterm tensor least squares problem
\begin{equation}\label{ttprob}
    \min_{\cc{X}} \| \cc{F} - \sum_{i =1}^{\ell} \cc{X} \times_1 A_1^{(i)} \times_2 A_2^{(i)} \cdots \times_d A_d^{(i)} \|_F,
\end{equation}
where $\cc{X}\in\bb{R}^{m_1 \times m_2 \times \cdots \times m_d}$, $\cc{F}\in\bb{R}^{n_1\times n_2 \times \cdots \times n_d}$ are tensors with $d$ dimensions (or modes), while $A_j^{(i)}\in\bb{R}^{n_j \times m_j}$ are tall matrices, for  $i=1,\ldots,\ell$ and $j=1,\ldots,d$. The term $\cc{F}$ is assumed to be in
low-rank Tucker format.

Matrix and tensor formulations of least squares problems have emerged in the recent literature as an alternative to classical vectorized forms in different data science problems, see, e.g., \cite{Bouseeetal.17},\cite{Brandoni.Simoncini.20},\cite{Chang.Wu.24},\cite{Kernfeldetal.15},\cite{doi:10.1137/110842570},\cite{algarte2025}. In particular,  the occurrence of a multiterm coefficient operator has been first proposed in \cite{Dantas2019},\cite{Dantasetal2017}. The numerical literature on the topic is very scarce, although the problem is challenging because of the absence of direct methods that can efficiently handle multiple addends in tensorial form, already for problems with small dimensions. The difficulty emerges as soon as the number of modes $d$ is equal to three or larger, giving rise to the so-called curse of dimensionality problem. Computationally, the solution of the linear algebra problem becomes intractable. Among the strategies available in the literature,
the LSQR method has been recently explored for $\ell=1$ (one addend in
(\ref{ttprob})) in \cite{Bentbibetal.22}, where Tucker and CP decompositions are used.

The problem of numerically solving tensor matrix equations with {\it square} coefficient matrices has attracted a lot of attention
in the past few decades, due to its growing occurrence in many scientific applications \cite{Hackbusch.Khoro.Tyrty.05},\cite{Khoromoskij.12},\cite{Khoromskij.18}.
However, also in the square case, the Tucker format is hard to handle without resorting
to a complete unfolding of the modes, so that
direct procedures
are mostly confined to low values  of $d$; see, e.g., \cite{Chen.Kressner.20},\cite{Simoncini.00}.
Except for special sparse tensor structures, see, e.g.,
\cite{Kressner.Tobler.10}, iterative methods suffer from the same memory problems.

With the classical Tucker format being too expensive, both in terms of computational costs and memory consumption, the tensor-train format (TT-format) has gained a lot of consideration, thanks to more affordable memory consumption for large $d$. This setting has been confirmed by our computational experience on \eqref{ttprob}, hence we will mainly focus on the TT-format for storing all quantities of interest, except for the coefficient matrices, which will be kept in their original form (the 
Kronecker product is never explicitly generated).
Tensor-train formulations have been successfully used in
the past decade in a variety of application problems, 
 directly approximating the given tensor operator, or
 for solving linear algebra problems in tensor format;
 see, e.g., \cite{doi:10.1137/110833142},\cite{doi:10.1137/140953289},
 and
 \cite{Buccietal.24} for a recent
 summary of contributions.
 The use of the TT-format is often
 accompanied by alternating least squares strategies, where
 a minimization with respect to single portions of the
 tensors is performed in turn. Here we take a different
 direction, by trying to stay as close as possible to the
 minimization of the whole problem, in a way that 
 deviations are only due to the need to maintain a memory saving
 TT-structure.

Our contribution is twofold. First, we provide the implementation of a tensor-train version of the LSQR algorithm, a popular method introduced by Paige and Saunders in \cite{paige1982lsqr} for solving the large and sparse least squares problem $\min_{{\bm x}} \|{\bm f}-{\bm A}{\bm x}\|$.
In the vector case, it is known that LSQR is preferable to the Conjugate Gradient (CG) on the normal equation because of stability issues, though the two methods are mathematically equivalent.
Thanks to the available matlab software package
TT-Toolbox \cite{software_tt}, the tensor implementation is very natural, and 
the tensor-train format can be easily maintained
throughout the iteration using specifically designed functions. 
We also propose the inclusion in our TT-LSQR method of a randomized strategy via
Johnson-Lindenstrauss transforms, to lighten the computational
costs when the problem becomes very large.
Our second contribution is the investigation of the problem
formulation (\ref{ttprob}) in addressing a recurrent problem
in information retrieval associated with term-document data:
the allocation of a new query document among
groups of documents, which were
already clustered using keyword terms as discriminant.
Our computational experience shows that the approach is
very promising and could be easily adapted to very large
datasets by exploiting the features of the
randomized projection.

An outline of the paper follows. 
In section \ref{sec:prel} we recall some tensor properties
and notation, mostly focusing on the
tensor-train (TT-) format.
Section~\ref{sec:tensor_lsqr} describes the tensor-oriented
implementation of the LSQR method, while section~\ref{sec:ttlsqr}
specializes to the TT-format. Section~\ref{sec:speed} and its subsections
briefly describe the procedures we have considered to enhance the 
computational performance of the developed method, that is, preconditioning and
randomized sketching.
 As a preliminary ``sanity check'' step, the new
algorithm and its features are explored on square problems in 
section~\ref{sec:sanity}.
In section~\ref{sec:IR} the problem of interest in the context of Information Retrieval is introduced. In section~\ref{sec:discr_expes} and its subsections
a large number of experiments on the allocation problem
with several benchmark datasets is reported.
In particular, in section~\ref{sec:expes_sketch} we describe how to exploit
sketching strategies to make the whole procedure computationally
more efficient, while the conclusions are drown in 
section~\ref{sec:conclusions}. The Appendix reports some
relevant aspects associated with the Tucker format, in case
that tensor formulation is of interest.

\section{Preliminaries}\label{sec:prel}
A tensor is a multidimensional array $\cA\in\RR^{n_1\times \dots \times n_d}$, denoted by capital italic letters to distinguish them from matrices. The tensor is said to be $d-$dimensional if it has $d$ dimensions (or \textit{modes}).
One important operation with tensors is the so-called $j$-mode product: given a tensor $\cA\in\RR^{n_1\times \dots\times n_d}$ and a matrix $U\in\RR^{m\times n_j}$, the $j$-mode product is denoted by
\[
\cB = \cA\times_j U\in\RR^{n_1\times\dots\times m\times\dots\times n_d},
\]
defined as follow
\[
\cB(i_1,\ldots,i_{j-1},k,i_{j+1},\ldots,i_d)=\sum_{i_n=1}^{n_j}\cA(i_1,\ldots,i_j,\ldots,i_d)U(k,i_j).
\]
The properties of the $j$-mode product are described in \cite{doi:10.1137/07070111X}. We recall the correspondence between the $j$-mode products and the Kronecker products. Given a tensor $\cX\in\RR^{n_1\times\dots\times n_d}$ and the matrices $A^{(j)}\in\RR^{m_j\times n_j}$, $j=1,\ldots,d$, we have
$$
\cY=\cX\times_1 A^{(1)}\times_2 A^{(2)} \dots \times_d A^{(d)}\quad\Leftrightarrow\quad
\by=(A^{(d)}\otimes \dots\otimes A^{(1)})\bx,
$$
where $\by$ and $\bx$ denote the vectorization of the tensors $\cY$ and $\cX$, respectively.

\begin{figure}[htb]\label{fig_tucker}
\centering
    \includegraphics[width=0.4\textwidth]{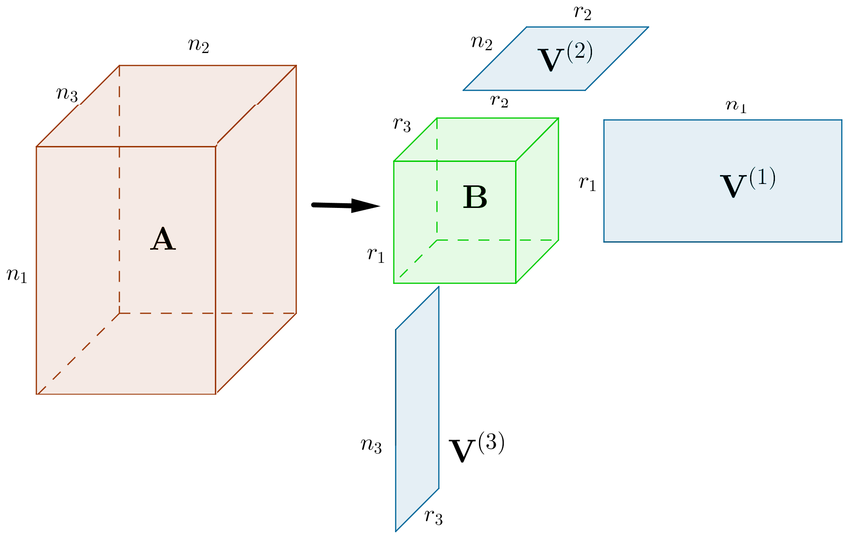}
    \caption{Tucker format.}
\end{figure}

Given the introduction of the $n$-mode product, the Tucker decomposition is a natural associated tensor
factorization.
\begin{definition}
Given a tensor $\cX\in\RR^{n_1\times \dots\times n_d}$, the Tucker decomposition of $\cX$ is defined as
\[
\cX\approx \cc{C}\times_1 A^{(1)}\times_2 A^{(2)}\dots \times_N A^{(d)},
\]
where $\cc{C}\in\RR^{m_1\times \dots\times m_d}$ is the core tensor and $A^{(j)}\in\RR^{n_j\times m_j}$.
\end{definition}
If the matrices $A^{(i)}$ of the Tucker decomposition are full column rank, the decomposition is said to be \textit{independent}; if the matrices $A^{(i)}$ have orthonormal columns, the decomposition is said to be \textit{orthonormal}. In Figure \ref{fig_tucker} an example of Tucker decomposition of a $3$-dimensional tensor is shown.

Alternatively, a very popular tensor representation is the Tensor-Train (TT) decomposition.
A tensor $\cc{X}$ is in TT-format if 
\begin{equation}\label{def_tt}
    \cc{X}(i_1,i_2,\ldots,i_d) = G_1(i_1,:)G_2(:,i_2,:)\cdots G_d(:,i_d),
\end{equation}
where each $G_k(:,i_k,:)\in\bb{R}^{r_{k-1}\times r_k}$, under the assumption $r_0=r_d=1$.

We can interpret each $G_k(:,i_k,:)$ as a matrix, and the whole of $G_k$ as a three-dimensional tensor, with dimensions $r_{k-1}\times n_k \times r_k$. Consequently, each element of $\cc{X}$ can be seen as
\begin{eqnarray}\nonumber
    \cc{X}(i_1, i_2,\ldots, i_d) &=& \sum_{\alpha_0,\alpha_1,\ldots,\alpha_d}G_1(\alpha_0, i_1, \alpha_1)G_2(\alpha_1, i_2, \alpha_2)\cdots G_d(\alpha_{d-1}, i_d, \alpha_d)\\ 
    &=& \sum_{\alpha_1,\ldots,\alpha_{d-1}}G_1( i_1, \alpha_1)G_2(\alpha_1, i_2, \alpha_2)\cdots G_d(\alpha_{d-1}, i_d).\\ \nonumber
\end{eqnarray}

 \begin{figure}[htb]
\centering
    \includegraphics[width=0.4\textwidth]{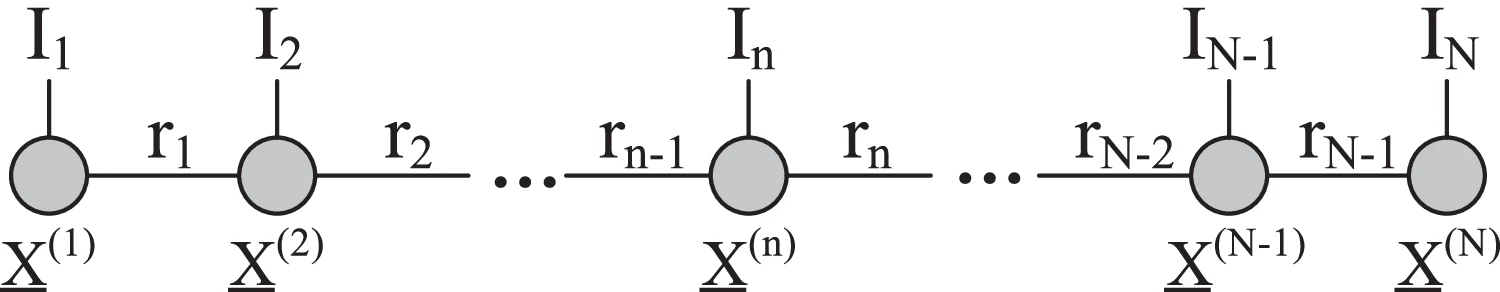}
    \caption{Tensor-Train format.}
\end{figure}

In the Tensor-Train format, a tensor is low-rank if each TT-rank is small, $r_i\ll n_i$. Classically, low rank is obtained through the \textit{rounding} operation, based on the (truncated) TT-SVD algorithm, introduced by Oseledets in \cite{Oseledets.2011}. Therefore, rounding has the role of what we call truncation function for other formats. We refer to section~\ref{sec:ttlsqr} for a few more details on the TT-SVD.



\section{Truncated Tensor-Oriented LSQR}\label{sec:tensor_lsqr}

In this section, we introduce the numerical method for solving the tensor least squares problem (\ref{ttprob}).

A straightforward way to solve (\ref{ttprob}), under the assumption $m = m_1 = \ldots = m_d$, employs vectorization through the Kronecker product, that is 
$$
{\bm A} =\sum_{i=1}^{\ell} A^{(i)}_d \otimes A^{(i)}_{d-1} \otimes \ldots \otimes A^{(i)}_1\in\bb{R}^{(n_1 n_2\cdots n_d) \times m^d}, 
\quad {\bm f}={\rm vec}(\cc{F})\in\bb{R}^{n_1 n_2 \cdots n_d}
$$ 
and ${\bm x}={\rm vec}(\cc{X})\in\bb{R}^{m^d}$,
and solves the vectorized least squares problem with
classical iterative LS methods. This procedure is inefficient due to the final dimensions of ${\bm A}$, which become infeasible even for a small initial problem size. 
Nonetheless, the Kronecker formulation is useful to analyze the spectral properties of the problem.

For $d=2$ the problem is reduced to a matrix least squares problem. For $\ell =2$ a truncated matrix-oriented LSQR algorithm has been recently studied in \cite{simoncini:hal-04437719}, where the implementation and properties of the truncated matrix-oriented LSQR are introduced to solve the generalized Sylvester least squares problem
\begin{equation}
    \min_{X\in\bb{R}^{m\times m}} \| F_1F_2^T - A_1^{(1)} X A_2^{(1)} - A_1^{(2)} X A_2^{(2)}\|_F,
\end{equation}
with $A_1^{(j)},(A_2^{(j)})^T\in\bb{R}^{n\times m}$ and $F_1,F_2\in\bb{R}^{n\times p}$, $p\ll \min\{n,m\}$. 
The method in \cite{simoncini:hal-04437719} is equipped with a computational efficient strategy to generate 
updated iterates of low rank. 
 Because of better numerical properties and a more efficient truncation function, truncated LSQR outperforms truncated CG on the normal equation, when solving matrix least squares problems.

We derive the truncated tensor-oriented LSQR keeping in mind the matrix-oriented procedure, although it is immediately clear that the level of complexity is much higher. The first step substitutes vectors with tensors,
\[
    {\bm x}^{(k+1)}={\bm x}^{(k)}+\frac{\phi_{k+1}}{\rho_{k+1}} {\bm d}^{(k+1)} \hspace{.3cm}\Longrightarrow\hspace{.3cm} \cc{X}^{(k+1)} = \cc{X}^{(k)} + \frac{\phi_{k+1}}{\rho_{k+1}} \cc{D}^{(k+1)}.
\]
With this preliminary step, the computed iterates are the same as those obtained in the vectorization case, but formally reshaped as tensors.
In fact, after the sum above the updated tensor needs to
 preserve the structure, while keeping memory
 requirements under control by ensuring that the rank remains small. Since a sum generally increases the rank, some truncation is necessary, like in the matrix case. 
 To adapt the aforementioned truncation strategy to the tensor case, it is necessary to take
 into account the tensor formats being used. We will mainly focus on the
 Tensor-Train (TT) format introduced in Section \ref{sec:prel}, keeping in mind that 
 the Tucker format has a different strategy for maintaining the low-rank tensor format of sums of tensors. 
 Summarizing, every time a tensor update is
 performed the sum requires to be computed
 within the chosen format, and then truncation
 is applied, again following the 
 format characterization. The same occurs
 when the coefficient operator and its 
 transpose is applied. 

 Preliminary computational experiments showed that the
Tucker formulation is far from being competitive for the class
of data we are interested in, both in terms of memory requirements and computational costs, already with
$d=3$ modes; this fact has been largely observed in the
recent literature. For this reason, in the
following we focus on the description of the tensor-train
properties, while the discussion of the Tucker format
is included for completeness in the Appendix.

\begin{algorithm}[htb]
\caption{TT-LSQR ALGORITHM}\label{tensor:lsqr}
\begin{algorithmic}
\REQUIRE $A^{(i)}_j$ with $j=1,\ldots,d$ and $i = 1,\ldots, \ell$, $\cF$, imax
\STATE $\beta_1 \cU_1=\cF$, $\alpha_1\cV_1=\cc{L}^T(\cU_1)$, $\cc{G}_1=\cV_1$, $\cc{X}_0=\cc{0}$, $\bar{\Phi}_1=\beta_1$, $\bar{\rho}_1=\alpha_1$
\WHILE{$i<{\rm imax}$}
\STATE $i=i+1$
\STATE $\beta_{i+1}\cU_{i+1}=\cc{L}(\cV_i)-\alpha_i\cU_i$
\STATE $\alpha_{i+1}\cV_{i+1}=\cc{L}^T(\cU_{i+1})-\beta_{i+1}\cV_{i}$
\STATE $\rho_i=(\bar{\rho}_i^2+\beta_{i+1}^2)^{1/2}$
\STATE $c_i=\bar{\rho}_i/\rho_i$, $\, \, s_i=\beta_{i+1}/\rho_i$
\STATE $\theta_{i+1}=s_i\alpha_{i+1}$,  $\qquad\bar{\rho}_{i+1}=-c_i\alpha_{i+1}$
\STATE $\Phi_i=c_i\bar{\Phi}_i$,
$\qquad\quad \bar{\Phi}_{i+1}=s_i\bar{\Phi}_i$
\STATE $\cc{X}_i=\cX_{i-1}+(\Phi_i/\rho_i)\cc{G}_i$
\STATE $\cc{G}_{i+1}=\cV_{i+1}-(\theta_{i+1}/\rho_i)\cc{G}_i$
\STATE \textbf{test for convergence}
\ENDWHILE
\end{algorithmic}
\end{algorithm}

\subsection{Tensor-train LSQR}\label{sec:ttlsqr}
The implementation of tensor-oriented LSQR in tensor train format, TT-LSQR in the sequel, is presented in Algorithm~\ref{tensor:lsqr}. In there, calligraphic Roman letters refer to tensors, while Greek letters refer to scalars (as an exception, $c_i, s_i$ also refer to scalars). The operations have to be adapted to the tensor setting:  section~\ref{sec:ttlsqr} discusses the tensor-train format, while the appendix reports the Tucker format. Here we just anticipate that thanks to 
particularly advanced software packages for tensor computations, the actual implementation does not substantially differ from 
Algorithm \ref{tensor:lsqr}, since all elementary operations 
involving tensors are overwritten, and thus
carried out according to their 
format\footnote{The code associated with Algorithm \ref{tensor:lsqr} will be made available in \\
{\tt https://github.com/Lorenzo-Piccinini/}
in the near future.}.Moreover, extra truncation can be applied, called ``rounding'' in the tensor-train format, to reduce the inner tensors size. We will
return to this issue later in this section.

In the TT-format the concept of low-rank TT-tensor is inherent in the format itself.
Theoretical results to estimate the error occurring when operating in the TT-format are available. In particular, if we assume that the given tensor is only approximately low-rank, i.e. the unfoldings satisfy
\begin{equation}\label{lowrank}
    A_k = R_k+E_k, \hspace{.2cm} {\rm rank}(R_k) = r_k, \hspace{.2cm} \|E_k\|_F=\epsilon_k, \hspace{.2cm}k=1,\ldots,d-1,
\end{equation}
then the error can be estimated as follows.
\begin{theorem}{\rm (\cite{OSELEDETS201070})}.
    Suppose that the unfoldings $A_k$ of the tensor $\cc{A}$ satisfy (\ref{lowrank}). The TT-SVD computes a tensor $\cc{B}$ in the TT-format with TT-ranks $r_k$ and
    \begin{equation}
        \|\cc{A}-\cc{B}\|_F^2\le \sum_{k=1}^{d-1}\epsilon_k^2 ,
    \end{equation}
    where $\epsilon_k$ is the distance (in the Frobenius norm) from $A_k$ to its best rank-$r_k$ approximation 
    \[
    \epsilon_k = \min_{{\rm rank}(B)\le r_k}\|A_k - B\|_F.
    \]
\end{theorem}
Usually, an SVD adapted to the TT context is used to approximate a full tensor with another one in TT-format (see below). In our case, the tensors are given in TT-format, but with sub-optimal TT-ranks due to the addition, which means that the tensor can be represented by another TT-tensor having lower TT-ranks. The TT-SVD algorithm can be adapted to these purposes. The complexity in this case is greatly reduced, while the error estimate does not change \cite{Oseledets.2011}. 
The procedure takes the name of what we have already called \textit{rounding}.

The considerations above are necessary to appreciate how operations such as sums and products are affected by rounding inside the LSQR algorithm. Indeed, whenever two TT-tensors $\cA,\cB\in\bb{R}^{n_1\times \cdots \times n_d}$ are added, the resulting tensor will have cores defined as
\[
C_1(i_1) = \begin{bmatrix} A_1(i_1) &  B_1(i_1))\\ \end{bmatrix}, \hspace{.3cm} C_d(i_d)=\begin{bmatrix}
    A_d(i_d) \\ B_d(i_d)\\ 
\end{bmatrix},
\]
and 
\[
C_k(i_k) = \begin{bmatrix}
    A_k(i_k) & 0 \\ 0 & B_k(i_k)\\ 
\end{bmatrix}, \hspace{.3cm} k = 2,\ldots, d-1.
\]
The sum tensor $\cc{C} = \cA + \cB$ belongs to $\bb{R}^{n_1\times \cdots \times n_d}$ and, if $r^A_i$ and $r^B_i$ are the TT-ranks of $\cA$ and $\cB$ respectively, the TT-ranks of $\cc{C}$ are $r^C_i = r^A_i + r^B_i$. Hence, whenever two or more tensors are added, the TT-ranks increase, and rounding becomes necessary to control the ranks. 
These truncations are taken care of in the Tensor-Train package \cite{software_tt}, by using a tensor-train
SVD.
%
This TT-SVD algorithm (\cite{Oseledets.2011}) is presented as an algorithm to compute the TT decomposition of a given tensor $\cA$. If the tensor is already in TT-format, its complexity is reduced. In particular, given a tensor $\cA\in\bb{R}^{n_1\times \cdots \times n_d}$ with suboptimal TT-ranks $r_k$, we want to estimate the reduced ranks $\bar{r}_k$ given an accuracy $\epsilon$ to preserve. We follow the
description in \cite{Oseledets.2011}: Recalling that
\[
\cA(i_1, i_2, \ldots, i_d) = G_1(i_1) G_2(i_2) \cdots G_d(i_d),
\]
we first proceed with a Right-to-Left orthogonalization of the cores, and then compress the TT-ranks through the SVD decomposition. To simplify the description, let us consider a $3-$dimensional tensor $G\in\bb{R}^{r_{k-1}\times n_k \times r_k}$, indexed as $G(\alpha_{k-1},i_k,\alpha_k)$.  The matrix $G(\alpha_{k-1},i_k \alpha_k)$ denotes its unfolding of dimensions $r_{k-1}\times n_k r_k$ and  $G(\alpha_{k-1}i_k, \alpha_k)$ its unfolding of dimensions $r_{k-1}n_k \times r_k$. Starting from the $d$th core up to the $2$nd core, at the $k$th iteration we compute the ${\rm QR}$ decomposition of the $k$th core $G_k \in \bb{R}^{r_{k-1}\times n_k \times r_k}$
\[
[G_k(\beta_{k-1},i_k \beta_k), R_k(\alpha_{k-1},\beta_{k-1})] = {\rm QR_{rows}}(G_k(\alpha_{k-1}, i_k \beta_k)),
\]
such that
\[
G_k(\alpha_{k-1}, i_k \beta_k) = R(\alpha_{k-1}, \beta_{k-1}) G_k(\beta_{k-1}, i_k \beta_k),
\]
where ${\rm QR_{rows}}$ denotes the reduced ${\rm QR}$ decomposition applied to the transpose of a matrix, leading to $G_k\in\bb{R}^{r_{k-1}\times n_k r_k}$ and $R_k\in\bb{R}^{r_{k-1}\times r_{k-1}}$. The $(k-1)$th TT-core is then updated as $G_{k-1} = G_{k-1}\times_3 R_k$.

The second step consists of the compression of the computed cores. Starting from the $1$st core up to the $(d-1)$th core, at the $k$th iteration the core $G_k$ is reduced via a singular value decomposition, that is
 $[G_k(\beta_{k-1}i_k, \gamma_k), \Lambda_k, V_k(\beta_k, \gamma_k)] = {\rm SVD}(G_k(\beta_{k-1} i_k, \beta_k)$,
where $G_k(\beta_{k-1}i_k, \gamma_k)\in\bb{R}^{\bar{r}_{k-1}n_k \times \bar{r}_k}$, 
$V_k(\beta_k, \gamma_k) \in \bb{R}^{r_k \times \bar{r}_k}$, $\bar{r}_k < r_k$ and
$\Lambda_k \in \bb{R}^{\bar{r}_k\times \bar{r}_k}$.
Then the next core is updated as $G_{k+1} = G_{k+1} \times_1 (V_k \Lambda_k)^T$. 

At the end of the TT-SVD algorithm the new TT-cores have dimension $\bar{r}_{k-1} \times n_k \times \bar{r}_k$, where $\bar{r}_k$ are the ranks obtained from the reduced ${\rm SVD}$. Further details can be found in \cite{Oseledets.2011}.

\section{Computational devices to improve performance}\label{sec:speed}
The vector LSQR algorithm may show slow convergence, when the normal equation has an ill-conditioned
coefficient matrix. Computational devices such as preconditioning
have been proposed to enhance the method performance.
In the following section we recall how to extend the preconditioning strategy to our
structured coefficient matrix. 
To lower the computational costs of the whole procedure, in section \ref{sec:sketching}
we also consider recently
developed sketching strategies, which seem to adapt well to the sum of Kronecker operators.

\subsection{Preconditioning for Tensor-oriented LSQR}\label{sec:precond}
The use of tensors makes each LSQR iteration significantly expensive, and the cost increases at each iteration until the maximum rank threshold is reached. After that
the cost per iteration remains approximately constant, though 
high, until termination is reached.
To limit computational costs we can apply 
preconditioning, so as to decrease the total number of iterations performed.

For the least squares problem, we consider right preconditioning \cite{BenziTumaPrec},\cite{HowellLsqrPrec}. Given a least squares problem, 
\[
\min_{\bx} \| {\bm f} - \bA\bx\|,
\]
with $\bA\in\bb{R}^{n\times m}$, $\bx\in\bb{R}^{m}$ and ${\bm f}\in\bb{R}^n$, we look for a matrix ${\bm M}\in\bb{R}^{m\times m}$ called \textit{preconditioner} such that the matrix $\bA{\bm M}^{-1}$ satisfies $\kappa(\bA {\bm M}^{-1}) \le \kappa(\bA)$, where $\kappa(\cdot)$ denotes the condition number. 
The problem is transformed as follows,
\[
\min_{\by} \| {\bm f} - \bA{\bm M}^{-1}\by \|, \hspace{.3cm} \by^* = {\bm M} \bx^*.
\]
The final solution is recovered as
$\bx^* = {\bm M}^{-1}\by^*$.
The matrix ${\bm M}$ is chosen in such a way that solving linear systems with ${\bm M}$ is sufficiently cheap.

In our tensor least squares problem the matrix $\bA$ is replaced by the tensor
${\cal \bA} = \sum_{i=1}^{\ell} {\cal \bA}^{(i)} = \sum_{i=1}^{\ell} {A}^{(i)}_d \otimes A^{(i)}_{d-1}\otimes \cdots \otimes A^{(i)}_1$, $\bx = {\rm vec}(\cX)$ and ${\bm f} = {\rm vec}(\cF)$. Depending
on the problem, the single term ${A}^{(i)}_j$ may
have few columns. Therefore, it is convenient to
choose as $\bm M$ the upper triangular factor of the ${\rm QR}$ decomposition of selected coefficient matrices $\bA^{(i)} = {\bm Q}^{(i)} {\bm R}^{(i)}$, that is
\[
A^{(i)}_j = Q^{(i)}_j R^{(i)}_j \hspace{.3cm}\Longrightarrow\hspace{.3cm} {\bm R}^{(i)} = R^{(i)}_d \otimes R^{(i)}_{d-1} \otimes \cdots \otimes R^{(i)}_1. 
\]
In practice, for a fixed mode, we compute the upper triangular factor for all $i=1,\ldots,\ell$, and for each mode $j=1,\ldots,d$ we choose the $R^{(i)}_j$ with the lowest condition number.
The preconditioner will be 
\[
{\bm M} = R^{(i_d)}_d\otimes R^{(i_{d-1})}_{d-1} \otimes \cdots \otimes R^{(i_1)}_1,
\]
where 
$i_j = {\rm arg}\min_{i=1,\ldots,\ell} \kappa(R^{(i)}_j)$. The strategy takes into account the fact that all matrices for each mode
have similar data, so that $R^{(i_j)}_j$ may
be effective on all matrices of the $j$th mode.

If the data are sparse, the preconditioned matrices ${\cal \bA}^{(i)} {\bm M}^{-1}$ do not need to be explicitly built. On the other hand,
if the number of dense columns of ${\cal \bA}^{(i)}$ is small, the preconditioned matrices are worth being computed once for all before starting the procedure.  

\subsection{Sketching for performance}\label{sec:sketching}
Although the use of tensors allows us
to save memory, the matlab implementation of TT-LSQR
yields unsatisfactory CPU time performance. 
Rank truncation alleviates this problem, however
speed remains an important bottleneck.
We addressed this issue by using
randomized linear algebra to solve overdetermined
least squares problems \cite{Rokhlin.Tygert.08}.
More precisely,
we introduce an
operator ${\cal S} \, : \, \RR^n \,\to\, \RR^s$
with
$s = 2 d \bar m$, where $d$ is the number of
modes in coefficient operator, and $\bar m$ is the total
number of columns in each mode; hence, $s$ is twice the
overall number of data columns involved in the computation. 
The operator ${\cal S}$ acts as a projection
onto a vector subspace of $\RR^s$, significantly
reducing the dimensions. 
This corresponds to replacing the problem 
$\min_{\bx} \|\bm f - {\mathcal {\bm A}} \bx\|$ with
$\min_{\bx} \|{\cal S} (\bm f - {\mathcal {\bm A}} \bx)\|$.

For the definition of $\cal S$ we consider the
popular Johnson-Lindenstrauss transform, giving
rise to the so-called sketching operator 
(\cite{woodruff2014sketching},\cite{martinsson_tropp_2020}).
This transform is given
by the composition of a diagonal matrix $D$ of $\pm 1$ with probability $1/2$,
then a Hadamard or Fourier-type transform $H$, and then a row sampling operator $J$,
giving
${\cal S}(A) = \frac 1 {\sqrt{mn}} J H D A$.

Letting $\bx_0$ be the solution to the sketched problem,
it can be shown \cite{martinsson_tropp_2020} that for a fixed $\epsilon>0$, bounds of the type
$$
\|\bm f - {\mathcal {\bm A}} \bx_0\| \le \frac{1+\epsilon}{1-\epsilon} 
\min_{\bx} \|\bm f - {\mathcal {\bm A}}  \bx\|. \quad\mbox{for}\quad
s\sim n \log n / \epsilon^2 .
$$
For our setting, the procedure corresponds to projecting the 
space of keywords (for the term-document datasets) or the space of
pixels (for the image dataset) into a significantly smaller vector subspace.
We construct a single sketching operator $\cal S$ that acts on the matrix
of each mode, for all terms, and we apply the operator so that
the coefficient matrices become $\widehat A_j^{(i)} = {\cal S}(A_j^{(i)})$.
The reduced problem to be solved thus reads as follows 
\begin{equation}\label{eqn:S_LS}
\min_{\bx} \|{\cal S} (\bm f) - (\sum_{i=1}^\ell {\widehat A}_d^{(i)}\otimes \cdots
\otimes {\widehat A}_1^{(i)}) \bx\| .
\end{equation}
Thanks to the choice $s\ll n$, the solution of (\ref{eqn:S_LS}) significantly reduces the computational cost of
the classification. In section \ref{sec:expes_sketch} we report on our
computational experience with this strategy.

\section{Computational experiments with discretized PDEs}\label{sec:sanity}
In this section we present a preliminary experimental analysis of the Truncated Tensor LSQR with Tensor-Train format. 

We consider the discretization by centered finite differences of the partial differential
equation
$$
-\Delta u + 2 \exp(1-x) u_x = f, \quad u=u(x,y,z), \quad  (x,y,z)\in [0,1]^3 ,
$$
equipped with uniform Dirichlet boundary conditions. Here $f$ is constant and equal to 1.
 Upon discretization, the problem can be written as
\begin{eqnarray}\label{eqn:pde_eq}
\left ( T\otimes I \otimes I +  I\otimes T \otimes I +  I\otimes I \otimes (T+B) \right ) \bx = \bm f ,
\end{eqnarray}
where $T\in\RR^{n\times n}$ is the usual tridiagonal matrix corresponding to the one-dimensional 
approximation of the negative second order derivative, and $B\in\RR^{n\times n}$ is the centered finite difference
approximation of $2 \exp(1-x) u_x$. We use this square problem as a ``sanity check" test, since the linear equation itself would not require a method for least squares problems to be solved.
 The fact that the equation residual goes to zero allows us to easily focus on different features of the method, such as the dependence on the truncation parameters. 

\begin{figure}[htb]
\centering
\includegraphics[width=0.49\textwidth]{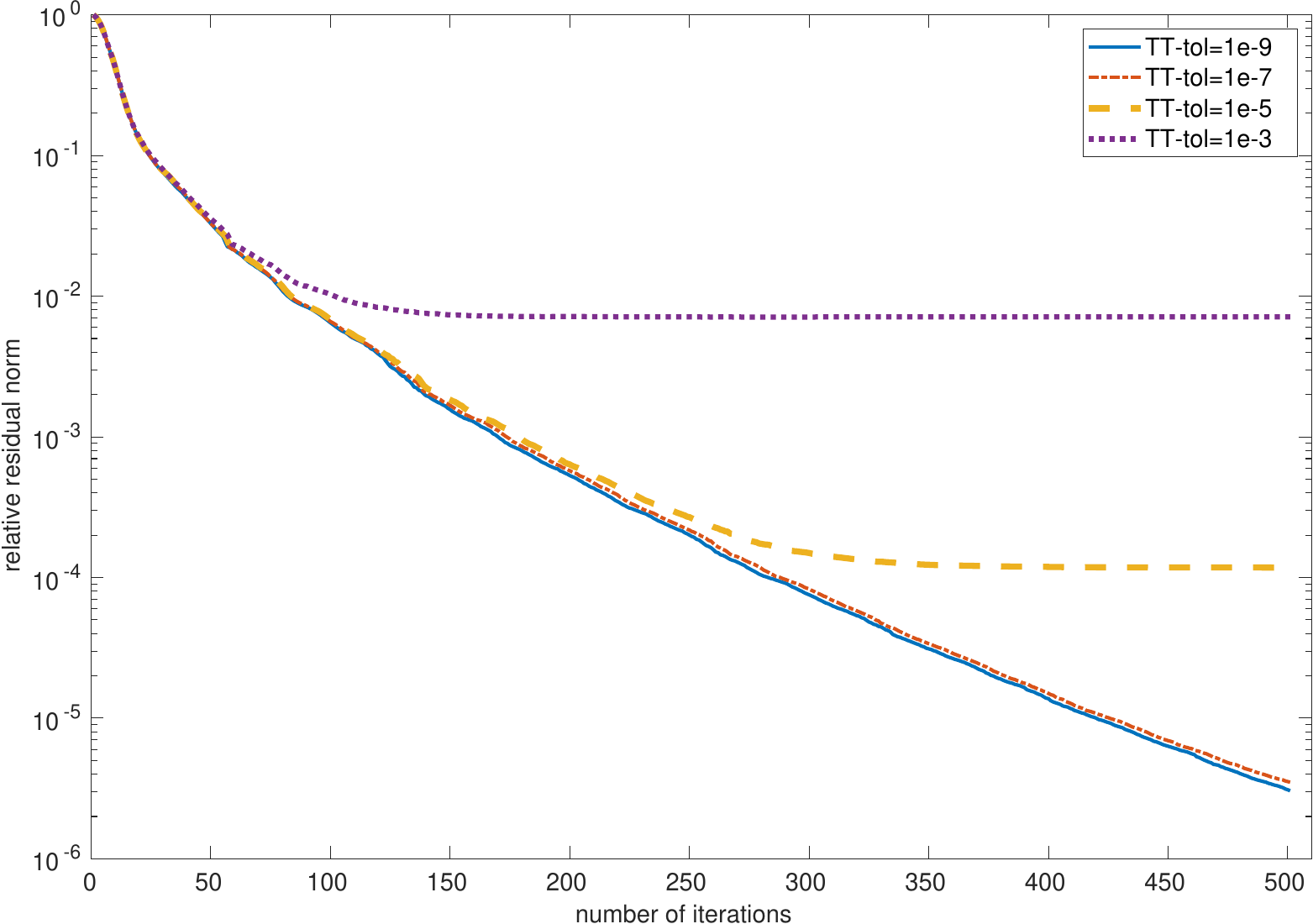}\,
\includegraphics[width=0.49\textwidth]{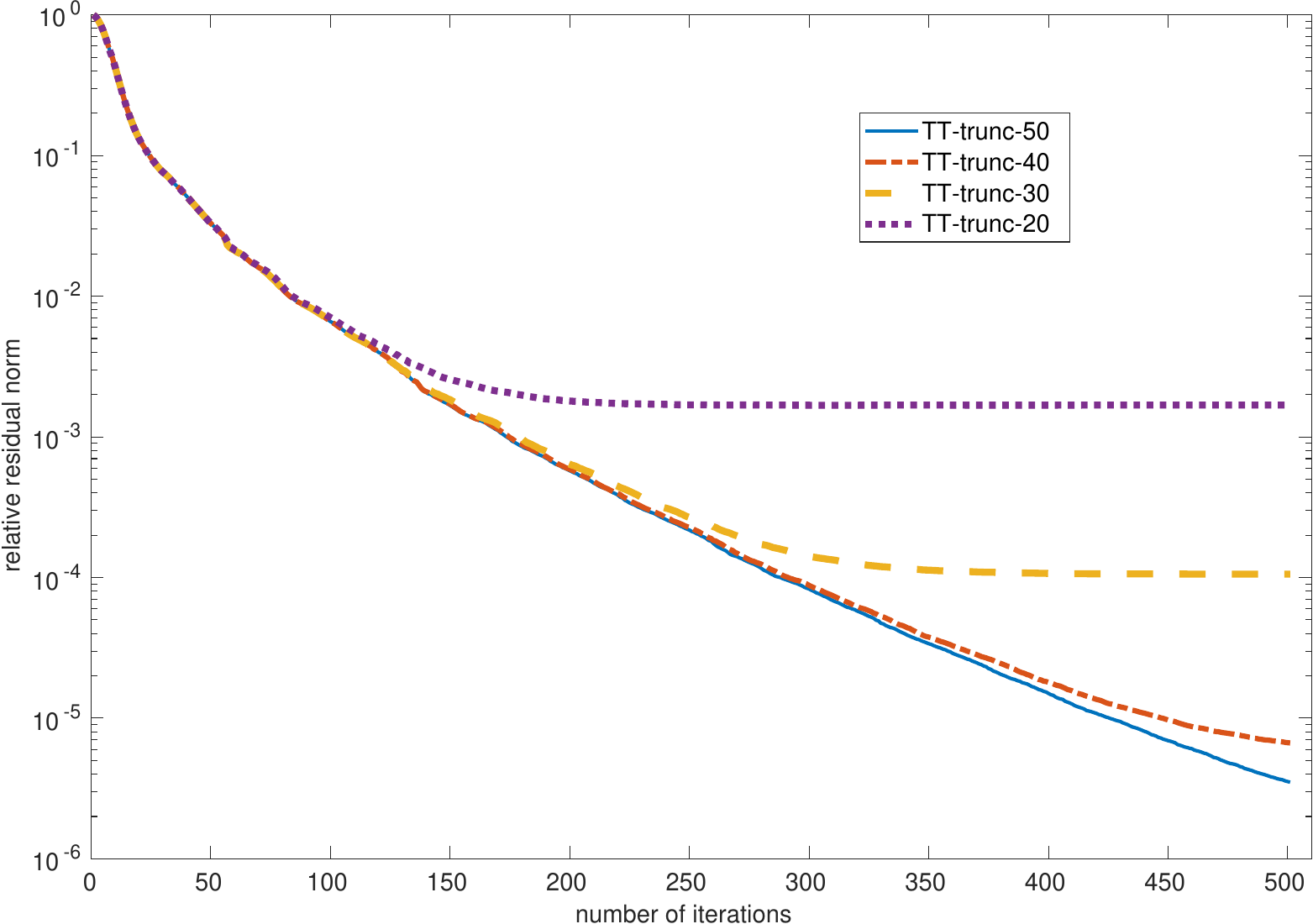}
\caption{Convergence history of TT-LSQR for $n=50$ for problem (\ref{eqn:pde_eq}). Left: Behavior as the TT truncation threshold varies. Right: Behavior as the TT truncation rank varies (truncation threshold $10^{-8}$).\label{fig:pde_param}}
\end{figure}

In the left plot of Figure \ref{fig:pde_param} we report the convergence history of our new method applied to (\ref{eqn:pde_eq}) for $n=50$, as the TT truncation parameter varies. As expected, stagnation occurs at a level that is related to the employed threshold. The results are in full agreement with what happens in the matrix case when the truncation threshold increases, see, e.g., \cite{Kressner.Tobler.11} for similar results using CG for symmetric and positive definite problems. 
The right plot of Figure \ref{fig:pde_param} shows the convergence 
behavior of the method with truncation threshold $10^{-9}$, as the maximum allowed rank of the TT-cores varies. Convergence is similarly affected. In general, truncation of information may be severe when applying the operator $\cal L$ with several terms, that is
for $\ell m$ large. 
Indeed, the TT-cores resulting after multiplications by $\cal L$ have size up to $\ell$ times the size of the argument tensor.  For this problem, the maximum rank equal to $n$ appears to be sufficient, and this is related to the fact that for this problem there is a single non-identity matrix for each mode, hence each TT-core behaves as if $\ell=1$. 

\begin{figure}[htb]
\centering
\includegraphics[width=0.49\textwidth]{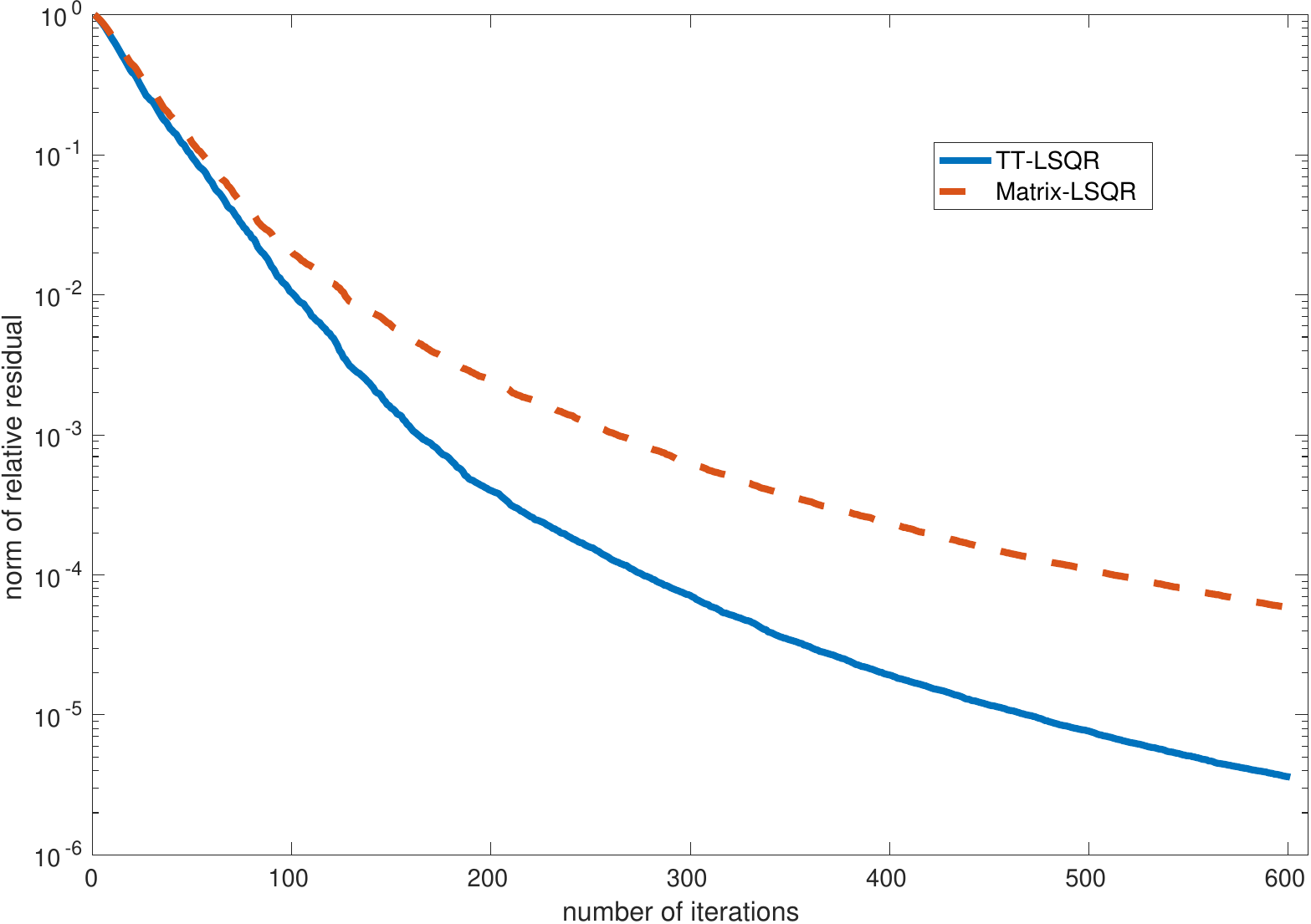}\,
\caption{Convergence history of TT-LSQR and MATRIX-LSQR
for $n=50$ for problems (\ref{eqn:pde_eq}) and (\ref{eqn:pde_eq2}), resp., for the discretization of (\ref{eqn:pde2}). Truncation threshold $10^{-9}$, rank truncation $50$.\label{fig:pde_conv}}
\end{figure}

The second experiment stems from the finite difference
discretization of a modification of the previous problem,
that is
\begin{equation}\label{eqn:pde2}
(a(x) u_x)_x + (a(y) u_y)_y +(a(z) u_z)_z + u_x = f,
\end{equation}
while the other settings are unchanged. Here
$a(x)=-\exp(-x)$, for $x\in (0,1)$, so that the
discretization of each second order derivative
evaluates $a(x)$ at the midpoint of each discretization
subinterval. The tensor equation takes again the form
(\ref{eqn:pde_eq}), where however $T$ and $B$ have
different meanings.
With this experiment we are interested in comparing the
performance of TT-LSQR with that of the matrix
version of the problem, which we call MATRIX-LSQR; for
the latter we use the matrix LSQR method 
developed in \cite{simoncini:hal-04437719}. 
The matrix-oriented formulation of the
problem can be obtained by collecting
one of the Kronecker products, giving
\begin{equation}\label{eqn:pde_eq2}
 TX  +  X\left ( T \otimes I +  I \otimes (T +
 B) \right )^T= \bm f_1 \bm f_2^T,
\end{equation}
with clear meaning for $\bm f_1, \bm f_2$; see, e.g., \cite{palitta2016matrix}.
Note that $X$ is a rectangular matrix of size
$n_3\times n_2n_1$. Given the structure of the
problem, in the matrix case truncation acts differently
than in the tensor case: all iterates have factors with either
$n_1$ rows or $n_2n_3$ columns each, hence memory requirements grow more significantly than in the tensor
case.
Figure~\ref{fig:pde_conv} reports the convergence history of
both methods for this problem. All matrices were
preconditioned as described in section~\ref{sec:precond}.
 TT rounding was performed
with tolerance $10^{-9}$ and rank truncation 50. The same
parameters were used for the truncation in MATRIX-LSQR. Truncation has a stronger negative
impact in the matrix formulation, due to the fact that more
information is lost after truncation.
In the plot, we can appreciate the better accuracy achieved
by the tensor method, in addition to the already mentioned
lower memory requirements.

\section{Information retrieval, discrimination and classification in data mining}\label{sec:IR}
Documents such as scientific articles, newspaper reports, or email messages can be categorized according to different factors. Among the most commonly used ones are terms and keywords used within the documents themselves.   In this framework, datasets are stored as
matrices, where the rows correspond to the term dictionary, while the
columns collect the distinct documents. The matrix entry $(i,j)$ is binary if the
term $i$ is contained at least once in the document $j$, or a natural number if the number of times the term $i$ occurs in the document $j$ is collected.

In {\it information retrieval}, given a term-document matrix, the problem consists of finding a good representation of a new document (query) in terms of those in the available dictionary. The problem
can be cast as a least squares problem \cite{EldenBook.19}. The largest components of the
solution vector give an indication of the documents closest to the test query. Very often, the dictionary documents have already
been clustered, so that they themselves are
already collected in groups having similar features, such as a similar topic. In this case, it may be relevant to {\it aggregate}
the new query as a member of one of the already separate clusters. In this case,
since the clusters have common features, 
it may be appropriate to talk about {\it discriminant} analysis, rather than
classification. Discriminant strategies have been used in multivariate statistics for a long
time \cite{Johnson.Wichern.07}, and they allow one to distinguish among elemental groups by using
features that are common within the group but different between groups.

On the other hand, 
classification is an important component of image processing, for which common/different features can be better identified than with text documents. Given a new image, corresponding
for instance to a person's face or to a specific class of objects (a car, a dog, etc),  the procedure determines to which person or class of objects the new image
belongs to. Each image is represented as a vector, a matrix or tensor of pixels.
In the case of vectors, images are collected as adjacent columns of a matrix. 
Different matrices represent different persons or objects.
Given a new image, comparing the corresponding vector with each different 
class using some distance measure allows one to classify the new image.
Our formulation fits naturally this representation, and it will be
used for one of the dataset in section~\ref{sec:discr_expes}.

A large number of classification methods, mostly for image processing,  have been studied. They exploit
different algorithmic devices, from the SVD to deep learning methods with
matrix and tensor techniques, see, e.g., \cite{Bouseeetal.17},\cite{Brandoni.Simoncini.20},\cite{doi:10.1137/110842570},\cite{Saputraetal.24} and
their references.
For both applications, discrimination and classification, 
 we are interested in exploring how the special additive structure of the
 tensor formulation can provide a surplus value to these methodologies.

\section{Discriminant analysis on clustered data}\label{sec:discr_expes}
In the following we focus on the allocation
of a new query, for already clustered data. 

An interesting aspect in using a sum of tensor formulation in the LS problem
for this type of data
analysis is the flexibility in
choosing the parameters $\ell$, $d$, and $m_i$ for $i=1,\ldots,d$. Fixed the number of modes $d$, if we want to use the same amount of memory, we can either choose to use more terms (larger $\ell$) with fewer columns (smaller $m_i$), or choose to have more columns in each term (larger $m_i$) but with fewer terms (smaller $\ell$). 
This has a multi-facet impact.
Assuming $m=m_1=\ldots=m_d$, the memory required in the tensor representation throughout the iteration depends on $m$. Focusing on the TT-format, we have $d$ core tensors with a memory requirement of the order of $\cc{O}(mr^2d)$, where $r$ is the largest TT-rank of the solution tensor. The overall amount of memory required will be of the order of $\cc{O}(n m d \ell + m r^2 d)$. If instead of $\ell$ terms we take only one term, to preserve the amount of information, the dimensions of $\cX$ should be increased to ${m} \times {m} \times \cdots \times {m}$, where $\bar{m}= m \ell$. Consequently, the new amount of memory required will be of the order of $\cc{O}(n {m} d + {m} r^2 d)$.
Using a bigger $m$ not only increases the memory requirements but, more importantly, also increases the computational cost of each operation. Each iteration requires multiple applications of the tensor operator $\cc{L}$, leading to several truncation steps, thus affecting the
final accuracy of the algorithm.
We refer to
Table~\ref{tab:m_ell} for an illustration of 
the effect of this selection on the method
performance.
Note that the selection of $d$ is not
left to the user, as it is mainly
problem dependent.
In all our experiments, TT rounding
was performed only using the tolerance, $10^{-4}$, and not the maximum rank dimension.

In the following we illustrate how the data are allocated to formulate the tensor LS problem, and the we briefly describe the main properties of the datasets. Several numerical experiments follow.

\subsection{Computational experiments: Preparing the setting}\label{sec:expes1}
In this set of experiments we consider term-document matrices that have already been grouped according to 
their common features. Assume we have $d=3$ groups, and
let ${\cal A}_1$, ${\cal A}_2$, ${\cal A}_3$
be the corresponding matrices, each having $\bar m$ columns. We then partition
the group matrices so that we have $\ell$ blocks of
size $m$
for each group, so that $m\ell=\bar m$, that is
\begin{equation}\label{eqn:calA}
{\cal A}_j=[A_j^{(1)}, \ldots,  A_j^{(\ell)}], \quad
A_j^{(i)}\in {\mathbb R}^{n\times m}, \quad
j=1,2,3.
\end{equation}
We are thus ready to define the tensor operator
$$
{\bm A}=\sum_{i=1}^\ell A_3^{(i)}\otimes A_2^{(i)}\otimes A_1^{(i)}.
$$
Given the query vector $f$, consisting of the recurring keywords of a new
document, we want to associate $f$ with the closest
group among the ones available, that is
${\cal A}_1, {\cal A}_2, {\cal A}_3$. To this end, we solve the
least squares problem
\[
\min_{{\bm x}}\| {\bm f} -{\bm A} {\bm x}\|_2
\]
with the new method,
where ${\bm f}=f\otimes f\otimes f$.

Once the structured least squares problem has been approximately solved, the classification of $f$ consists of identifying the group of documents closest to
$f$ by using the solution $\bx$. To this end, we investigate two criteria that exploit our new formulation.


\begin{description}
\item \underline{Criterion 1}\footnote{ with $\stackrel{j}{f}$ we mean that
$f$ is positioned in either of the three 
Kronecker terms.}:
$$
\hat j = {\rm arg}\max_{j=1, 2, 3} |(1\otimes  \stackrel{j}{f} \otimes 1)
\left (\sum_{i=1}^\ell A_3^{(i)} \otimes A_2^{(i)}\otimes A_1^{(i)} \right ) {\bm x}|,
\hskip 0.2in (C1).
$$
\item \underline{Criterion 2}:
\begin{equation}\label{dist:svd}
   \hat j = {\rm arg}\max_{j=1, \ldots, 3}
   \|f^T U^{(j)}\|, \hskip 0.2in (C2)
\end{equation}
where $U^{(j)}$ contains the left singular vectors of the unfolding in the mode $j$ of
$(\sum_{i=1}^\ell A_3^{(i)} \otimes A_2^{(i)}\otimes A_1^{(i)} ) \widetilde{\bm x}$,
and $\widetilde{\bm x}$ is an approximation to the TT solution $\bm x$ as described next. In the following wet ${\bm x}={\rm vec}({\cal X})$.
We first compute the TT-cores of the solution tensor, denoted $G^{(i)}\in\bb{R}^{r_{i-1}\times n_i \times r_i}$, then we compute their range basis using the ${\rm SVD}$. Since all the TT-cores, except for the first and the last one, are $3-$dimensional tensors, the ${\rm SVD}$ is performed with respect to the second mode, which coincides with the mode of the tensor we are working with, 
\begin{equation}
    [U^{(1)},\sim] = {\rm svd}(G^{(1)}), \hspace{.2cm} [U^{(i)},\sim] = {\rm svd}(\hat{G}^{(i)}), \hspace{.2cm} [U^{(d)},\sim] = {\rm svd}(G^{(d)^T}),
\end{equation}
where $\hat{G}^{(i)}\in\bb{R}^{n_i \times r_{i-1}r_i}$ denotes the unfolding along the second mode of the TT-core $G^{(i)}\in\bb{R}^{r_{i-1}\times n_i \times r_i}$.
This procedure becomes expensive for a larger number of modes (i.e., high $d$) and for large TT-ranks. Instead of using the full solution tensor $\cX$, the ${\rm HOSVD }$ of $\cX$ is computed and just the first column of each $U^{(i)}$ is stored. 
This means that $\cX(i_1,\ldots,i_d) \approx 
\widetilde \cX(i_1,\ldots,i_d):= U^{(1)}(i_1,:) \cdots U^{(k)}(:,i_k,:)\cdots U^{(d)}(:,i_d)$, using the Tensor-Train notation. Then we apply the operator $\cc{L}$ to $\widetilde \cX$, which will have TT-ranks equal to $\ell$, the number of terms. The TT-cores of the tensor $\cc{L}(\widetilde\cX)$ precisely correspond to the tensors $\widetilde{U}^{(i)}\in\bb{R}^{\ell \times n_i \times \ell}$, used to classify the image $f$ in the criterion (C2).
\end{description}

\vskip 0.1in
Criterion 2 is more expensive than Criterion 1,
as it requires the SVD computation of the unfolded terms. If no truncation in $\cX$ were performed,
we would expect (C2) to be  more reliable
as a discriminant criterion, as it may be viewed
as a tensor generalization of the usual subspace
projection criterion for data stored as matrices (see Criterion 4 below).
On the other hand, the truncation may significantly
weaken this strategy, with respect to (C1).
Indeed, we noticed in our experiments that (C2)
is not always more successful than (C1) in
recognizing clustered data. For this reason, one could consider a combination of the two criteria to refine the classification. A more detailed analysis
is left for future investigations.

We compare the quality of the success of the allocation with
classical procedures in data mining that do not involve solving a tensor LS problem: a test for the allocation of queries in information retrieval, somewhat similar to simple query matching \cite[section 3.1.2]{Berry.Browne.99}, 
and a test using range-oriented recognition \cite[section 11.2]{EldenBook.19}.
\vskip 0.1in
\begin{description}

    \item \underline{Criterion 3}:
Partition $\bm w = {\rm arg}\min_{\bm w} \|\bm f-[{\cal A}_1, {\cal A}_2, {\cal A}_3] 
{\bm w}\|_2$
as ${\bm w}=[w_1;w_2;w_3]$, according to the partition of the coefficient matrix. 
For this matrix setting, and exploiting that the coefficient matrix
has been normalized, we consider
the following elementary query matching criterion:

\vskip 0.1in
\begin{quotation}
Doc $f$ belongs to Category $\widehat j\,$ if  $\,\,\widehat j= \arg\max_j \|w_j\|_2$ \hskip 0.1in (C3).
\end{quotation}
\vskip 0.1in
 \item \underline{Criterion 4}: Let ${\cal U}_i$, $i=1, \ldots, 3$ be the
 matrix whose $10$ orthonormal columns span the leading portion of 
 {\rm range}$({\cal A}_i)$. The matrix ${\cal U}_i$ is obtained by
 a truncated SVD of ${\cal A}_i$ to the first 10 singular triplets \cite[Section 11.2]{EldenBook.19}.
 Then
 \vskip 0.1in
 \begin{center}
Doc $f$ belongs to Category $\widehat j\,$ if  $\,\,\widehat j= \arg\min_j \|f - {\cal U}_j{\cal U}_j^T f\|$ \hskip 0.1in (C4).
\end{center}
\vskip 0.1in
 
\end{description}
\vskip 0.1in
Criterion 3 requires solving a large matrix least squares problem,
while Criterion 4 requires computing $d$ truncated SVD of the
corrisponding matrices.

\begin{figure}[hbt]
\centering
\includegraphics[width=1.6in,height=1.6in]{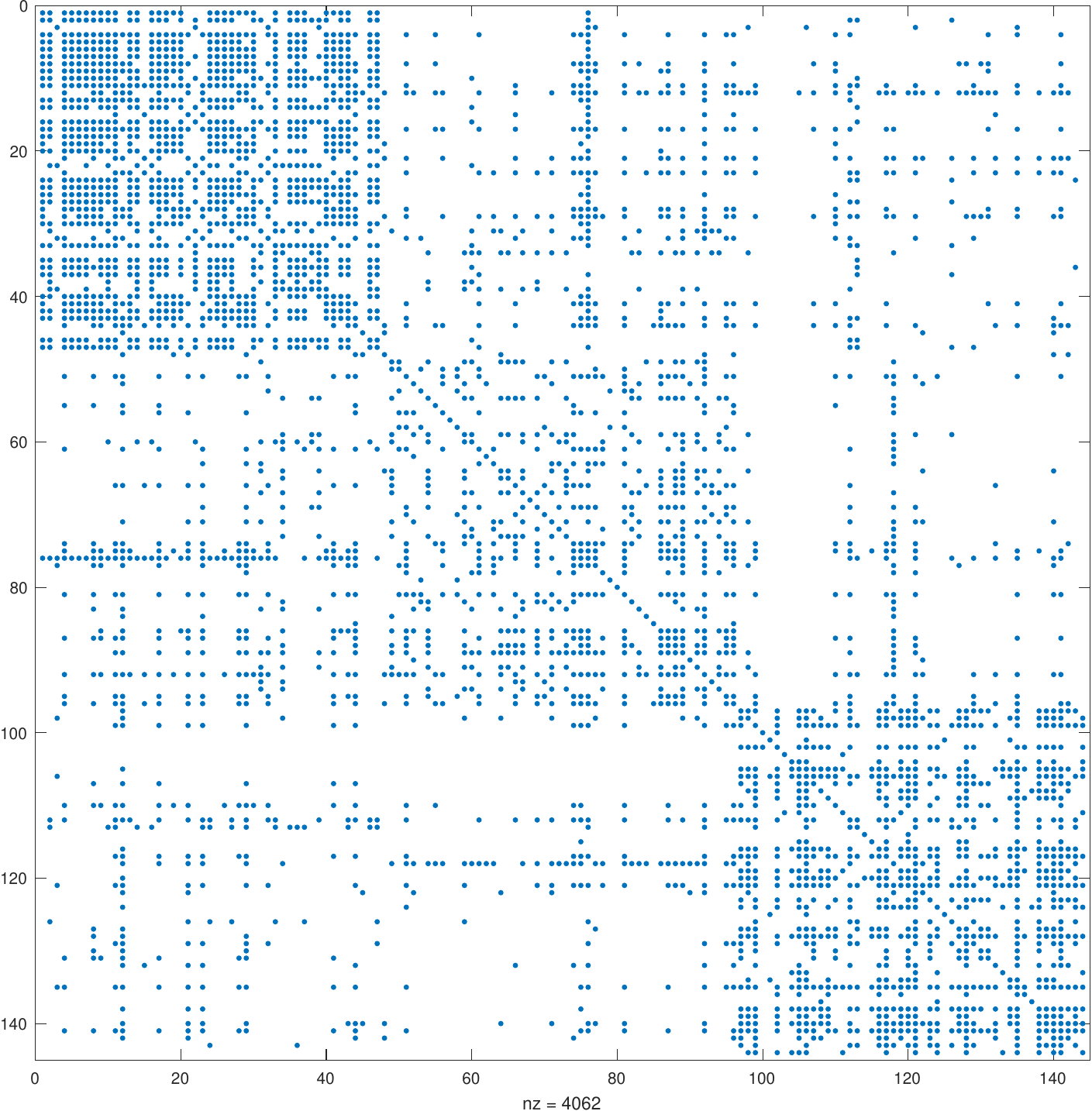} \,
\includegraphics[width=1.6in,height=1.6in]{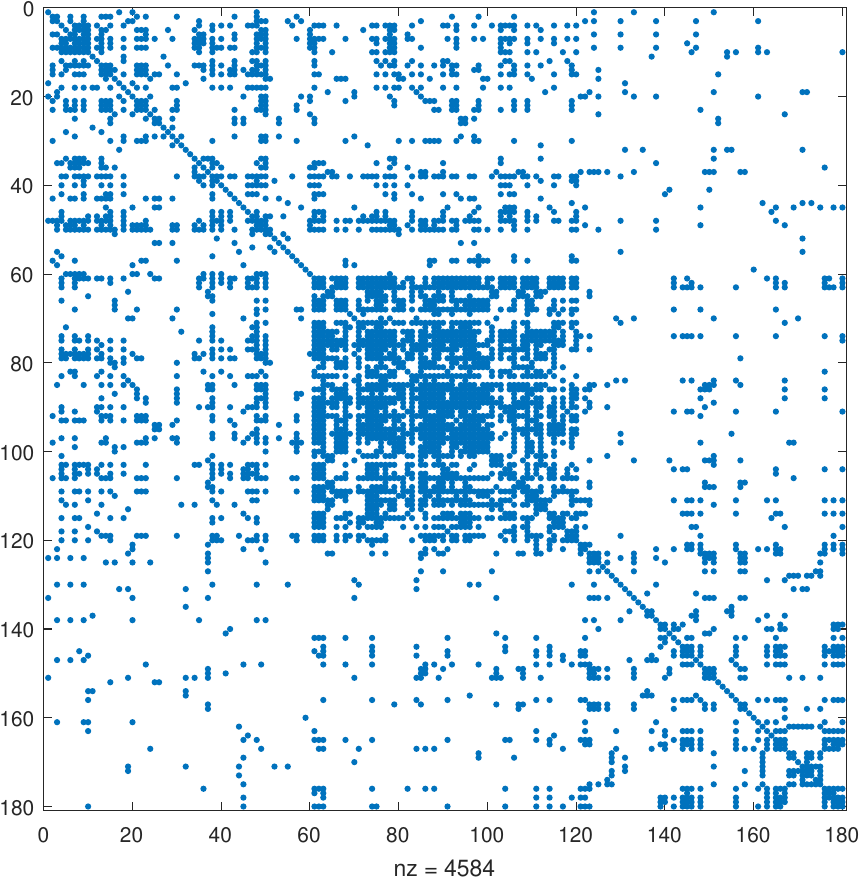}\,
\includegraphics[width=1.6in,height=1.6in]{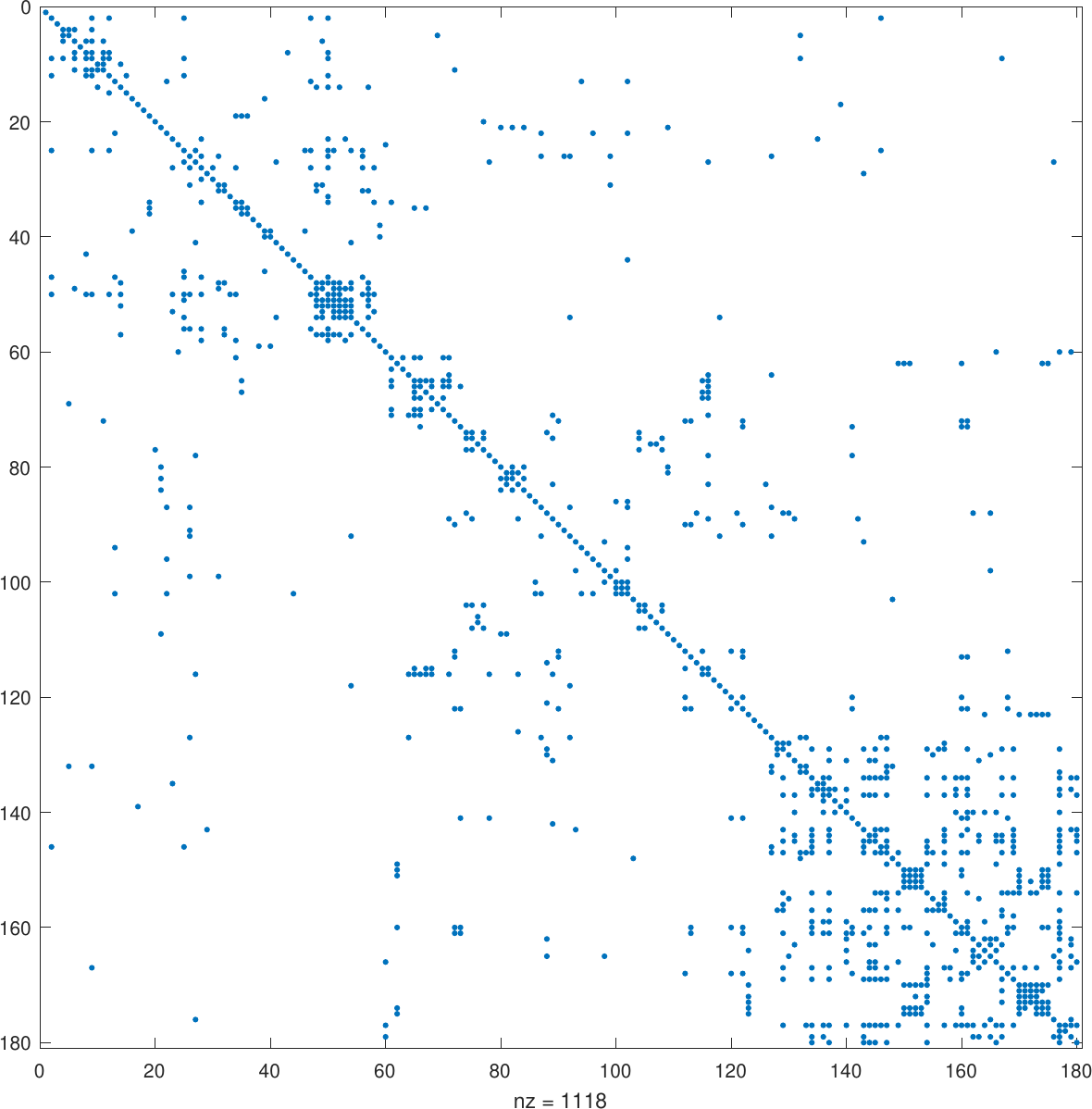} 
   \caption{Sparsity pattern of correlation matrices of
   ${\cal A}$ in (\ref{eqn:calA}) with $m=60$ for {\tt Reuters}(left), {\tt Cranfield}(middle), {\tt Medline} (right); reported are elements above 0.15 in absolute value.\label{fig:spydata}}
\end{figure}

\vskip 0.1in

 For our experimental investigation, we consider three term-document datasets available in popular repositories\footnote{These datasets
are available at {\tt https://ir.dcs.gla.ac.uk/resources/test\_collections/}}:

\vskip 0.1in
\begin{itemize}
\item {\tt Reuters21578}. This 
is a widely used collection of newswire articles
from the Reuters financial newswire service, and it is an essential benchmark for text categorization research, with 8293 documents and 18933 keywords, giving a very sparse matrix with a total of 389455 nonzero entries. The documents are already assembled and indexed with categories, so that 65
clusters are already identified. In the following
we will work with the three
largest groups;

\item {\tt Cranfield}. This is a term-document matrix of scientific contexts on physics topics with 1,400 articles
and a list of 4563 keywords, for a total of 81201 nonzero entries in the data matrix;

\item {\tt Medline}. This is a term-document matrix  of medical contents with 1,033  articles
 and 5735 keywords, for a total of 51174 nonzero entries in
the data matrix;

\end{itemize}
\vskip 0.1in 

\begin{figure}[hbt]
\centering
\includegraphics[width=2.3in,height=1.8in]{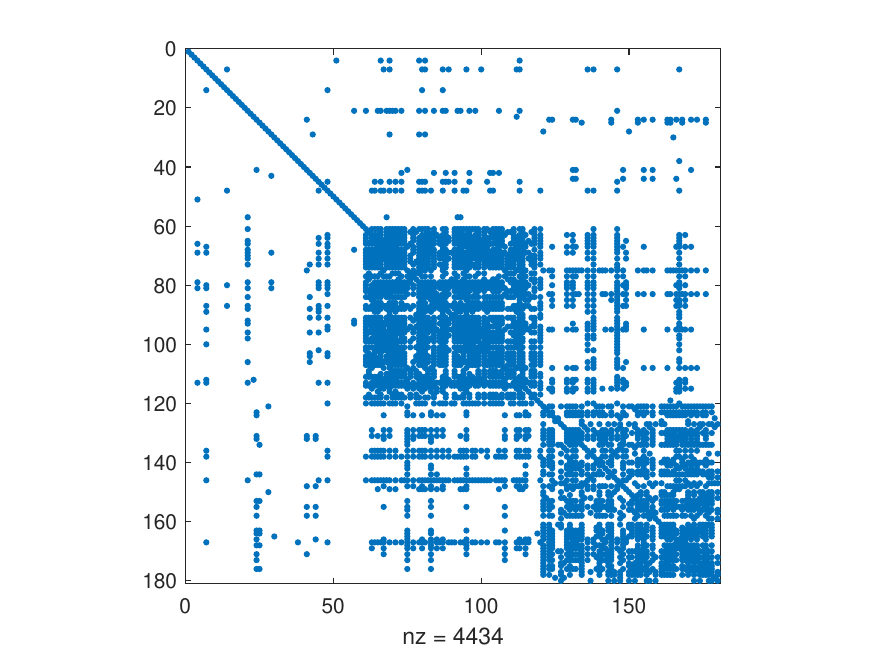} \,
   \caption{Sparsity pattern of correlation matrix of
   ${\cal A}$ in (\ref{eqn:calA})
   for {\tt FashionMNIST}; with $m=60$  reported are elements above 0.7 in absolute value.\label{fig:spydata1}}
\end{figure}

Except for {\tt Reuters21578}, the term-document
datasets need to be divided into clusters, as described above.
To this end, the K-means algorithm (available in Matlab) with three clusters was
employed  for each data set.  By taking the
computed clusters as categories, we have selected the first
$\bar m$ columns of each cluster as training sets, and
the subsequent 20 columns as test set.  Since the clustering is
empirical, we do not expect exceptional success rates,
although the results appear to be quite satisfactory in most cases.

We also consider a popular dataset of images:

\vskip 0.1in 
\begin{itemize}
\item {\tt Fashion\_MNIST}\footnote{The dataset is available at {\tt https://github.com/zalandoresearch/fashion-mnist}.}. This is a widely used dataset containing $70000$ images of fashion products from Zalando. The images, stored as pixel vectors of length $784=28^2$, 
are divided into Training and Testing sets. As a preprocessing, we create two matrices of dimension $784\times n_{\tt train}$ and $784\times n_{\tt test}$ each, with $n_{\tt train} = 60000$ and $n_{\tt test} = 10000$. The products are already classified in $10$ categories each representing a different type of clothing. In our experiments we consider only the three groups of footwear: sandals, sneakers and boots, in this order.

\end{itemize}
\vskip 0.1in

\begin{table}
\centering
\begin{tabular}{|c|c|rr|r|r|}
\hline
 &  &  \multicolumn{2}{c|}{Tensor LS} & Matrix LS  & Proj norm\\
Category & $\bar m$ &   Crit (C1) & Crit (C2) & method (C3) &  (C4)\\
\hline
1 & 36 &  95 &100 & 100 & 100 \\
  & 48 &  95 &100 & 100 & 100\\
  & 60 & 90  &100 & 100 & 90 \\
\hline
2 & 36 &95 & 85 &25 & 85 \\
  & 48 & 85 & 95  & 10 & 85 \\
  & 60 & 90 & 95  &  0 & 95  \\
  \hline
3 & 36 &95 & 95 &90 & 95\\
  & 48 &95 & 95  &  85 & 85 \\
  & 60 &95  & 95   & 85 & 95   \\
  \hline
    \end{tabular}
    \caption{Dataset {\tt Reuters}. Percentage of success of classification for 20 test documents for different groups.  $\ell=6$, $d=3$, $n=18933$ for varying $\bar m=\ell m$.\label{tab:reuters}}
\end{table}

Figure~\ref{fig:spydata} shows the sparsity pattern of the correlation matrix associated with the
term-document matrices,
while  Figure~\ref{fig:spydata} shows the sparsity pattern relative to
{\tt FashionMNIST}.
We can observe that some of the groups  are
strongly correlated, shown by a dense diagonal block (e.g., the second group of
{\tt Cranfield}). However, a good classification
also depends on the sparsity of the off-diagonal
blocks, which corresponds to low correlation between groups (e.g., the third group in {\tt Cranfield}).
We also notice the high sparsity of the
{\tt Medline} correlation matrix, particularly for all entries of the second group, predicting a possible low quality of a classification procedure.
Similarly for the first group (sandals) of 
{\tt Fashion\_MNIST}.

\begin{table}
\centering
\begin{tabular}{|c|c|rr|r|r|}
\hline
 &  &  \multicolumn{2}{c|}{Tensor LS} & Matrix LS  & Proj norm\\
Category & $\bar m$ &   Crit (C1) & Crit (C2) & method (C3) &  (C4)\\
\hline
1 & 36 &  55 & 65&10 & 75 \\
  & 48 & 80 & 80 & 25 & 95 \\
  & 60 & 65 & 95 & 40 & 85 \\
\hline
2 & 36 &100 & 100 &100 &95\\
  & 48 &100 & 90 & 100 &90 \\
  & 60 &95 &  85 & 100 & 80  \\
  \hline
3 & 36 &85 & 95 & 10 &95 \\
  & 48 & 95 & 85 & 25 & 80 \\
  & 60 & 75  & 85 &10 & 70  \\
  \hline
    \end{tabular}
    \caption{Dataset {\tt Cranfield}. Percentage of success of classification for 20 test documents for different groups.  $\ell=6$, $d=3$, $n=4563$ for varying $\bar m=\ell m$.\label{tab:cran}}
\end{table}

\begin{table}
\centering
\begin{tabular}{|c|c|rr|r|r|}
\hline
 &  &  \multicolumn{2}{c|}{Tensor LS} & Matrix LS  & Proj norm\\
Category & $\bar m$ &   Crit (C1) & Crit (C2) & method (C3) &  (C4)\\
\hline
1 & 36 &  100 &95 & 90 &95 \\
  & 48 & 95 & 85 & 75 & 90 \\
  & 60 & 90 & 90 & 80 & 85 \\
\hline
2 & 36 &  80 &70 & 10 & 60 \\
  & 48 & 55 & 55 & 10 & 25 \\
  & 60 & 60 & 60 & 5 & 55 \\
  \hline
3 & 36 & 95 & 95 & 95 & 95 \\
  & 48 & 95 & 95 & 100 & 95 \\
  & 60 & 100 & 100 & 100 & 100 \\
  \hline
    \end{tabular}
    \caption{Dataset {\tt Medline}. Percentage of success of classification for 20 test documents for different groups.  $\ell=6$, $d=3$, $n=5736$ for varying $\bar m=\ell m$.\label{tab:med}}
\end{table}

\begin{table}
\centering
\begin{tabular}{|c|c|rr|r|r|}
\hline
 &  &  \multicolumn{2}{c|}{Tensor LS} & Matrix LS  & Proj norm\\
Category & $\bar m$ &   Crit (C1) & Crit (C2) & method (C3) &  (C4)\\
\hline
1 & 36 & 35 & 70 & 5 & 60 \\
 & 48 & 65 & 75 & 0 & 75 \\
 & 60 & 75 & 75 & 0 & 60 \\
\hline
2 & 36 & 100 & 80 & 90 & 80 \\
 & 48 & 100 & 85 & 95 & 100 \\
 & 60 & 95 & 95 & 100 & 95 \\
  \hline
3 & 36 & 50 & 100 & 100 & 100 \\
 & 48 & 30 & 95 & 70 & 95 \\
 & 60 & 30 & 95 & 55 & 90 \\
  \hline
    \end{tabular}
    \caption{Dataset {\tt Fashion\_MNIST}. Percentage of success of classification for 40 test documents for different groups.  $\ell=6$, $d=3$, $n=784$ for different values of $\bar m=\ell m$.\label{tab:fashmnist}}
\end{table}

\subsection{Computational experiments: Testing the collocation strategy}\label{sec:expes2}
The results of our analysis are reported in Tables \ref{tab:reuters}-\ref{tab:fashmnist},
where the percentage of success in correctly
classifying 20 test documents is reported.
A fixed number of 10 iterations of TT-LSQR was used, as a very accurate solution 
is not the goal for this application. The TT
truncation parameter was set to $10^{-4}$.
 The tables show the performance as the total
 number $\bar m=\ell m$ of 
 documents in the test group increases, with
 fixed $\ell=6$.
Since the test sample is small, the difference in success between 85\% and 90\% is given by a single extra failure, therefore a difference in 5\% is considered not sufficiently large to make a case.

Table \ref{tab:reuters} shows that all test samples
are quite well identified, and that tensor-based (C2) is the best criterion. On the other hand, the simple matrix LS criterion (C3) is quite unreliable, giving exceptionally bad rates for group 2. This
performance is common to other datasets.
The results for the other term-document tests lead to similar conclusions, with a general preference
of criterion (C2) over (C1) for the tensor strategy.
We recall that except for {\tt Reuters}, the term-document groups were constructed by a clustering pre-processing, so that classification errors are more likely.

Table \ref{tab:fashmnist} shows that for the  image dataset, the first group is the hardest to recognize. This behavior is justified by the low correlation of the first group, as shown in Figure \ref{fig:spydata1}. The low correlation in the first group can be attributed to the significant differences among the images. Nonethess, the tensor strategy performs well, equipped with criterion (C2).
Finally, we observe that in most cases
the classification performance does not seem to over-depend on the value of
$\bar m =\ell m$ for any of the strategy,
while $\bar m=60$ is often (but not always)
the most favourable. As already mentioned,
some of the differences may depend on the
strength of the category cluster.

We would like to linger over the role of $\ell$
and $m$, that is the data distribution among the available
terms in the coefficient sum of tensors. For the application problem considered,
 the $m$ columns of each mode can be differently distributed among the
 terms: $m$ columns for each of the $\ell$ terms, so that $m \cdot \ell = \bar m$.
 For fixed $\bar m=60$,
 Table~\ref{tab:m_ell} illustrates the significant benefit, in terms of CPU time, of keeping $m$ low,
 and use more terms in the sum, without
 sacrifying classification success. The displayed results refer to the dataset {\tt Medline},
 after solving problem (\ref{ttprob}) to classify 20 documents from the first group, and reach a tolerance for the
 residual norm of the normal equation less than $10^{-4}$. A maximum of 200 iterations was allowed.
 In addition to significantly lower CPU times, using a small value of $m$,
 that is thin matrices $A_j^{(i)}$, allows one to work with lower rank tensors 
 throughout the TT-LSQR computation, making the iterations lighter. It is also
 interesting that convergence to the required
 residual norm tolerance is much faster for
 small $m$. These two properties lead to  dramatically lower
 CPU times for smaller $m$.



\begin{table}
\centering
\begin{tabular}{|c|c|c|c|c|}
\hline
 $\ell$ & $m$& \# iter (avg) & CPU time (secs, avg) & (C2)  \\
\hline
\,\,3 & 20 & \!\!*200 & 1277.3 &  60\\
\,\,6 & 10 & \,\,35 &\,\,257.3 &  80 \\
12 & \,\,5 & \,\,18 & \,\,190.4 &  85 \\
\hline
\end{tabular}
\\
\caption{Dataset {\tt Medline}. Performance of TT-LSQR {over 20 test documents} as $\ell$ and $m$ vary and are such
that $\bar m=\ell m  =60$. Stopping criterion on the normal equation relative
residual with tol=$10^{-4}$, and 200 maximum iterations. (*) Only 10 test documents are considered, due to excessive computational costs.\label{tab:m_ell}}
\end{table}

\subsection{Computational experiments with sketching}\label{sec:expes_sketch}
 In Table~\ref{tab:sketch} we report the
performance of the sketched versus the unsketched tensor least squares
problem, using a fixed number of iterations of TT-LSQR, equal to 10. Two of
the datasets are considered.
The percentage of successful classification is reported. For the first category,
CPU times are also included; since times are very similar for the other categories,
they are not reported in the table. We notice the drastic reduction of CPU time
for the sketched problem, with in most cases an acceptable decrease in performance.

In spite of the convenient bound in (\ref{eqn:S_LS}), the sketched solution may
significantly differ from the minimizer of the original problem.
It was recognized in \cite{doi:10.1137/23M1551973}, \cite{Rokhlin.Tygert.08}
that the sketched solution $\bx_0$ could be used as starting guess for an
iterative procedure in terms of the original least squares problem.
We have exploited this strategy in our setting as follows: we ran 30 iterations
of the sketched problem at low cost, and then "regularized" the obtained
solution by 2 iterations of TT-LSQR by solving the problem
\begin{equation}\label{eqn:twopass_LS}
\min_{\bm z} \|(\bm f-{\mathcal {\bm A}}  \bx_0) - {\mathcal {\bm A}}  \bm z\|.
\end{equation}
Our final solution will then be $\bx_k=\bx_0 + \bm z_k$. We explicitly remark
that the idea of using $\bx_0$ as starting guess suggested in \cite{doi:10.1137/23M1551973}, \cite{Rokhlin.Tygert.08} requires
that the sketched problem be solved exactly. To meet this constraint, we used 30 iterations of TT-LSQR on the sketched problem, and indeed running only 10 iterations was not sufficient to ignite the subsequent two iterations for  (\ref{eqn:twopass_LS}).
We remark that without truncation
the starting residual $\bm f -{\mathcal {\bm A}}  \bx_0$ has TT-rank larger than
$\bm f$, therefore the two iterations of TT-LSQR may be expensive.
The computational results of this strategy are reported in the two rightmost columns of Table~\ref{tab:sketch} for each
dataset.
We see that the percentage of success 
often pays off the extra cost, compared with the purely sketched
strategy. We believe this combined strategy deserved
a deeper analysis.


\begin{table}
\centering\footnotesize
\begin{tabular}{|c|c|l|l|l|l|l|l|}
\hline
Cate- & $\bar m$ &\multicolumn{3}{c|}{\sc medline}  &\multicolumn{3}{c|}{\sc Cranfield} \\
gory &  & {Tensor}& {Tensor} & Tensor w/& {Tensor}& {Tensor} & Tensor w/\\
 &  &   Sketched  & & StartSketch&   Sketched  & & StartSketch\\
\hline
1 & 36 &  90 (0.7) &95 (\,\,38.9) & 85 (13.9)& 35 (0.7)& 65 (27.2)& 45 (11.2)\\
  & 48 & 80 (1.8) & 85 (\,\,76.5) &85 (26.2)& 85 (2.0) &80 (58.9)& 70 (21.2)\\
  & 60 & 85 (4.3) & 90 (121.9)& 75 (50.9)& 90 (4.4) &95 (96.5) &  65 (33.3)\\
\hline
2 & 36 &  45  &70  & 85  &100 & 100& 100 \\
  & 48 & 45  & 55  & 55 & 95 & 90 & 100 \\
  & 60 & 65  & 60  & 75  &85  & 85& 95 \\
  \hline
3 & 36 & 85   & 95  & 90& 60 & 95  & 90 \\
  & 48 & 85 & 95  &90 & 95 &85  & 70 \\
  & 60 & 90  & 100  &100  &90  & 85 & 75 \\
  \hline
    \end{tabular}
    \caption{Percentage of success of classification for 20 test documents for different categories, using
    criterion (C2). Comparison
    of Sketched TT-LSQR (10 its), TT-LSQR (10 its) and  TT-LSQR (2 its) with the sketched solution (30 its)
    as starting guess. $\ell=6$, $d=3$, $n=5736$ for varying $\bar m$. Sketching parameter 
    $s=2 d \bar m$. In parenthesis is the average CPU time in secs - it is similar
    for all other categories. \label{tab:sketch}}
\end{table}

\section{Conclusions}\label{sec:conclusions}
We have proposed a tensorized version of the iterative
method LSQR for solving tensor-oriented multiterm least
squares. The implementation fully exploits available software
packages for the Tensor-Train representation of the
data, considerably lightening technical 
crucial steps such as rank truncation.
We have focused our analysis on the use of tensor multiterm
least squares in the context of allocating new document queries,
and on the advantage of having a solver that can directly
deal with this least squares formulation.

The obtained classification results are very promising and compare quite
well with standard strategies, for instance in the case of the less exercised
term-document datasets.

\section{Acknowledgments}
We would like to thank Davide Palitta for explanations on
\cite{Buccietal.24}.
 Both authors are members of the INdAM Research
Group GNCS. 

The work of VS
was partially supported by the European Union - NextGenerationEU under the National Recovery and Resilience Plan (PNRR) - Mission 4 Education and research
- Component 2 From research to business - Investment 1.1 Notice Prin 2022 - DD N. 104 of 2/2/2022,
entitled “Low-rank Structures and Numerical Methods in Matrix and Tensor Computations and their
Application”, code 20227PCCKZ – CUP J53D23003620006.
The work of LP is funded by 
the European Union under the National Recovery and Resilience Plan (PNRR) - Mission 4 - Component 2 Investment 1.4 ``Strengthening research structures and creation of ``National R\&D Champions'' on some Key Enabling Technologies'' DD N. 3138 of 12/16/2021 rectified with DD N. 3175 of 18/12/2021, code CN00000013 - CUP J33C22001170001.

\section*{Appendix}
{\bf Tucker-LSQR}. 
Given a tensor $\cc{X}\in\bb{R}^{n_1\times \cdots \times n_d}$  in Tucker format, that is
\begin{equation}
    \cc{X} = \cc{C} \times_1 X_1 \times_2 X_2 \cdots \times_d X_d,
\end{equation}
with $\cc{C}\in\bb{R}^{r_1 \times\ldots \times r_d}$ and $X_i\in\bb{R}^{n_i \times r_i}$, an exact representation of
$\cc{X}$ is obtained for $r_i = n_i$. 
On the other hand, a Tucker tensor is low-rank if $r_i\ll n_i$ for all $i=1,\ldots,d$. When implementing the tensor LSQR with the
Tucker format, the truncation step compresses the dimensions of the core tensor of the given tensor through the (truncated) Higher-Order SVD (HOSVD) developed in \cite{doi:10.1137/S0895479896305696}.

\begin{figure}[htb]\label{fig:blockdiag}
\centering
\includegraphics[width=0.5\textwidth]{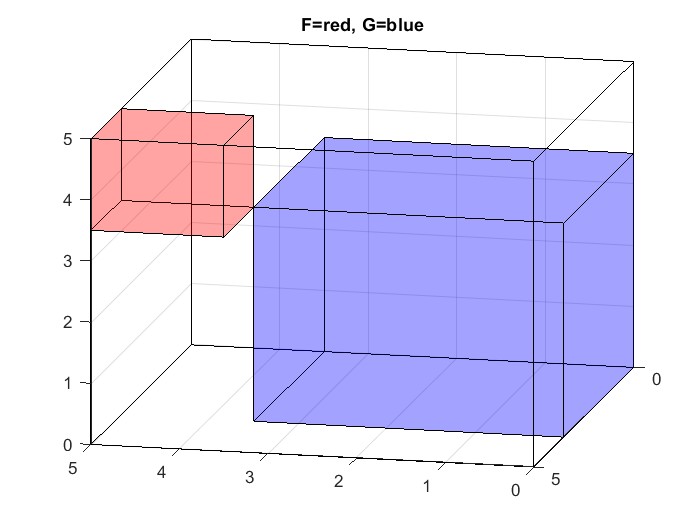}
\caption{Tensor core of the sum of two Tucker tensors,
$\cS={\rm blkdiag}(\cc{H},\cG)$.}
\end{figure}

The Tucker structure can be deployed to sum up two tensors without explicitly building them. Let $\cA=\cc{H}\times_1 U_{1}\times_2 U_{2}\dots\times_d U_{d}$ and $\cB=\cG\times_1 V_{1}\times_2 V_{2}\dots \times_d V_{d}$ have consistent dimensions, 
with $\cc{H}\in\RR^{r_1\times r_2\times \dots\times r_d}$ and $\cG\in\RR^{p_1\times p_2\times \dots\times p_d}$.
The sum tensor is $\cA + \cB = \cS \times_1 W_1 \times_2 W_2 \cdots \times_d W_d$, where
\[
W_{i}=[U_{i}, V_{i}],
\]
and, setting $q_1\times q_2\times \dots \times q_d$ where $q_i = r_i + p_i$,
\begin{eqnarray*}
\cS(1:r_1,1:r_2,\ldots,1:r_d)&=&\cc{H} ,\\ 
\cS(r_1+1:q_1,r_2+1:q_2,\ldots,r_d+1:q_d)&=&\cG ,
\end{eqnarray*}
which corresponds to the usual notation 
$\cS={\rm blkdiag}(\cc{H},\cG)$ (here in 3-mode tensor form); 
see Figure \ref{fig:blockdiag}.

If no further action is taken, the core tensor dimensions of all constructed tensors increase within each iteration. Truncation is therefore necessary to control memory requirements. The matrix truncation strategy developed in
\cite{simoncini:hal-04437719}
can be generalized to this setting.

Another advantage of the Tucker format follows from the
following theorem.

\begin{theorem}{\rm (\cite{jiang2017tensor})}\label{th_tuck}
For a given tensor $\cX\in\RR^{n_1 \times \dots\times n_d}$ with independent Tucker decomposition $\cX= \cc{C}\times_1 A^{(1)}\times_2 A^{(2)},\ldots\times_d A^{(d)}$, we have that
\begin{equation}
    \alpha\|\cc{C}\|_F\le\|\cX\|_F\le\beta\|\cc{C}\|_F,
\end{equation}
where $\beta=\|A^{(1)}\|_2\|A^{(2)}\|_2\dots\|A^{(d)}\|_2$ and $\alpha=\beta/(\Pi_{i=1}^{d}\kappa(A^{(i)})$ with $\kappa(A^{(i)})$ being the condition number of the matrix $A^{(i)}$.
\end{theorem}
The result above shows that if for all $i$ the matrix $A^{(i)}$ has orthonormal columns, the equality $\|\cX\|_F=\|\cc{C}\|_F$ holds. More generally,
the norm of the whole tensor can be monitored by that of the core. 

\bibliographystyle{siam}
\bibliography{references}

\end{document}